\documentclass{article}
\usepackage{amsmath,amsthm,amsfonts,amscd,eucal}

\begin{document}
\title{ Operator algebras and vertex operator algebras}
\author{Sebastiano Carpi
\\
Dipartimento di Economia, \\
Universit\`a ``G. D'Annunzio" di Chieti-Pescara,\\
Viale Pindaro, 42, I-65127, Pescara, Italy \\
E-mail: s.carpi@unich.it \\ 
}
\date{}
\maketitle

\begin{abstract} In two-dimensional conformal field theory (CFT) the building blocks are given by chiral CFTs, i.e.~CFTs on the unit circle (compactified light-ray). They are generated by quantum fields depending on one light-ray coordinate only. There are two mathematical formulations of chiral CFT, the one based on vertex operator algebras (VOAs) and the one based on conformal nets. We describe some recent results which, for first time, gives a general construction of conformal nets from (unitary) VOAs. 
\smallskip

\end{abstract}

\section{Introduction}

The study of CFT in two space-time dimensions has found applications to different areas of physics and mathematics such as string theory, quantum gravity, critical phenomena, infinite dimensional Lie algebras, number theory, finite simple groups, 3-manifold invariants, the theory of subfactors and noncommutative geometry \cite{Gannon,Witten,CHL}. Chiral CFTs, i.e. CFTs on the circle, are the building blocks of CFT. We have two different mathematical formulations of chiral CFT: {\it vertex operator algebras} (VOAs) and {\it conformal nets}. 

The notion of VOA first appears in the work of Frenkel, Lepowsky and Meurman \cite{FLM} as an important special case of the related notion of vertex algebra introduced by Borcherds \cite{Borcherds}. The main motivation was of mathematical nature, namely the construction of the moonshine module $V^\natural$ for the Monster group $\mathbb{M}$, the largest among the sporadic finite simple groups, in relation to the so called monstrous moonshine conjectures of Conway and Norton which relate the Monster with various special modular functions \cite{Gannon}.

A VOA (over $\mathbb{C}$) is a complex vector space $V$ together with a linear map (the \textit{state field correspondence}) $V\ni a \mapsto Y(a,z)$ satisfying certain physically motivated assumptions \cite{FLM,FHL,Kac}. 
The {\it vertex operators} $Y(a,z)$ should be considered as quantum fields on the unit circle $S^{1}:=\{z\in \mathbb{C}: |z|=1\}$ and the axioms for VOAs can be naturally viewed as an algebraic version of the Wightman axioms \cite{SW}, see Chapter 1 in \cite{Kac}. 

Conformal nets are the chiral CFT analogue of the Haag-Kastler nets in {\it algebraic quantum field theory} (AQFT), the approach 
to quantum field theory based on the theory of operator algebras on Hilbert spaces, namely C*-algebras and von Neumann 
algebras \cite{Haag,Kaw2015}.  Although the algebraic structures plays a very important role in the theory of operator algebras the functional analysis aspects are crucial so that the conformal net approach can be considered mainly analytic in contrast to the essentially purely algebraic nature of the VOA approach. One of the most relevant aspects of the theory of conformal nets is the beautiful and powerful interplay with some central topics in the theory of operator algebras and in particular with the theory of subfactors initiated by Jones \cite{Jones} which has been explicitly related to AQFT by the index-statistics theorem discovered by Longo \cite{Longo}. 

A conformal net $\mathcal{A}$ on $S^1$ is a map $I\mapsto \mathcal{A}(I)$ from the family of open proper intervals of $S^1$ 
into the family of von Neumann algebras acting on a fixed Hilbert space $\mathcal{H}$ (the \textit{vacuum Hilbert space}) satisfying certain physically motivated assumptions \cite{Kaw2015}. Recall that a von Neumann algebra is an algebra of bounded (linear) operators acting on a Hilbert space which contains the identity operator and it is closed under taking adjoints and weak limits.

Despite their significant mathematical differences the VOA approach and the conformal net approach exhibit their common ``physical root" through many structural similarities. Moreover, various interesting chiral CFT models can be considered from both point of view with similar outputs. For these and other reasons the two approach are generally believed to be essentially equivalent  but, until recently, no general result was known in this direction. 

In this contribution we report on some of the main results obtained in a joint work with Kawahigashi, Longo and Weiner which gives for the first time a general construction of conformal nets from (unitary) VOAs satisfying some natural extra assumptions, see \cite{CKLW} for more details and for the proofs.

\section{Conformal nets and VOAs}
In this section we recall the definitions of conformal net and of unitary VOA. 
Let $\mathcal{I}$ be the family of open non-dense non-empty intervals of $S^1$ and let $\mathrm{Diff}^+(S^1)$ be the infinite dimensional Lie group 
of orientation preserving smooth diffeomorphisms of $S^1$.  For any $I\in \mathcal{I}$, we consider the subgroup 
$\mathrm{Diff}(I) \subset \mathrm{Diff}^+(S^1)$ defined by
$$\mathrm{Diff}(I) :=\{\gamma \in \mathrm{Diff}^+(S^1): \gamma(z)=z\;\textrm{for all}\; z\in I'\}.$$
A (local irreducible) conformal net $\mathcal{A}$ is an inclusion preserving map $I\mapsto \mathcal{A}(I)$ from $\mathcal{I}$ into the family of von Neumann algebras acting on a fixed Hilbert space $\mathcal{H}$ satisfying the following assumptions

\begin{itemize}
\item[(i)] {\it Locality}: if $I_1 \cap I_2 =\emptyset$ then $[\mathcal{A}(I_1),\mathcal{A}(I_2)] = 0$.
\item[(ii)] {\it Conformal covariance}: there is a strongly-continuous projective unitary representation $U$ of $\mathrm{Diff}^+(S^1)$ on 
$\mathcal{H}$ such that $U(\gamma)\mathcal{A}(I)U(\gamma)^*=\mathcal{A}(\gamma I)$ for all $I\in \mathcal{I}$ and all 
$\gamma \in \mathrm{Diff}^+(S^1)$ and such that $U(\gamma) \in \mathcal{A}(I)$ for all $\gamma \in \mathrm{Diff}(I)$. 
\item[(iii)] {\it Positivity of the energy}: $U$ is a positive-energy representation, i.e.~the self-adjoint generator $L_0$ of the one-parameter rotation subgroup of $U$ (the {\it conformal Hamiltonian}) has non-negative spectrum. 
\item[(iv)] {\it Vacuum}: $\mathrm{Ker}(L_0) =\mathbb{C}\Omega$, where $\Omega \in \mathcal{H}$ (the {\it vacuum vector}) is a unit vector which is cyclic for the von Neumann algebra $\bigvee_{I \in \mathcal{I}}\mathcal{A}(I)$ generated by all the local von Neumann algebras $\mathcal{A}(I)$, 
$I \in \mathcal{I}$, i.e.~$\bigvee_{I \in \mathcal{I}}\mathcal{A}(I)\Omega$ is a dense subspace of $\mathcal{H}$.
\end{itemize} 

\bigskip

We now come to VOAs. In order to simplify the exposition in this contribution we will consider only the so called VOAs of {\it CFT type}. Let $V$ be a vector space. A formal series $a(z) =\sum_{n\in \mathbb{Z}}a_{(n)}z^{-n-1}$ with coefficients $a_{(n)} \in \mathrm{End}(V)$, $n\in \mathbb{Z}$, is called a {\it field} on $V$ if for every $b \in V$ we have $a_{(n)}b=0$ for $n$ sufficiently large. 

A vertex operator algebra (of CFT type) is a ${\mathbb{Z}}$-graded complex vector space 
$V=\bigoplus_{n\in \mathbb{Z}}V_n$, with $V_n$ finite dimensional for all $n \in \mathbb{Z}$, equipped with a linear map 
$a \mapsto Y(a,z) =\sum_{n\in \mathbb{Z}}a_{(n)}z^{-n-1}$ from $V$ into the space of fields on $V$ satisfying the following assumptions

\begin{itemize}
\item[(i)] {\it Locality}: if $a,b \in V$ then $(z-w)^N[Y(a,z),Y(b,w)] = 0$, for all sufficiently large positive integers $N$.
\item[(ii)] {\it Conformal covariance}: there is a vector $\nu \in V$ (the {\it conformal vector}) such that the coefficients of the corresponding
vertex operators $Y(\nu,z)=\sum_{n\in \mathbb{Z}} L_nz^{-n-2}$ satisfy the Virasoro algebra commutation relations
$$[L_n,L_m] = (n-m)L_{n+m} + \frac{c}{12}(n^3-n)\delta_{n,-m} 1_V$$
with central charge $c \in \mathbb{C}$. Moreover,  $V_n = \mathrm{Ker}(L_0-n1_V)$ for all $n\in \mathbb{Z}$ and 
$[L_{-1},Y(a,z)]=\frac{d}{dz}Y(a,z)$ for all $a\in V$. 
\item[(iii)] {\it Positivity of the energy}: the representation of the Virasoro algebra associated with the conformal vector $\nu$ is a positive-energy representation, i.e.~$L_0$ has non-negative eigenvalues so that $V_n=0$ for all $n<0$. 
\item[(iv)] {\it Vacuum}: $\mathrm{Ker}(L_0) =\mathbb{C}\Omega$, where $\Omega \in V$ (the {\it vacuum vector}) is a unit vector which satisfies 
$Y(a,z)\Omega|_{z=0} =a$ for all $a\in V$. 
\end{itemize}
We now discuss unitarity. Let $(\cdot|\cdot)$ be a (positive definite) scalar product on $V$ which is {\it normalized},  i.e.~$(\Omega|\Omega)=1$.   
We assume that the pair $\big(V,(\cdot|\cdot)\big)$ has {\it unitary Virasoro symmetry}, i.e.~that the Virasoro algebra acts unitarily on $V$.
 Then, we can define the adjoint vertex operators by 
$Y(a,z)^+=\sum_{n\in \mathbb{Z}}a_{(n)}^+z^{n+1}$, $a\in V$, where, for any $n$, $a_{(n)}^+$ is defined by $(a_{(n)}^+b|c)=(b|a_{(n)}c)$, 
$b,c \in V$.  Now, let $V$ be a VOA with a normalized scalar product  $(\cdot|\cdot)$. The pair $\big(V,(\cdot|\cdot)\big)$ is said to be a unitary VOA if it has unitary Virasoro symmetry and if for every $a\in V$ the adjoint vertex operator $Y(a,z)^+$ is mutually local with all vertex operators $Y(b,z)$, $b\in V$, i.e.~if for every pair $a,b \in V$ we have $(z-w)^N[Y(a,z)^+,Y(b,z)]=0$ for all sufficiently large positive integers $N$.  
An equivalent definition of unitarity can be given through the notion of invariant bilinear form for VOAs defined in \cite{FHL}. 
The latter definition has been also considered in \cite{DL}. The equivalence is a consequence of a VOA version of the 
PCT theorem proved in \cite{CKLW}. 

\section{From VOAs to conformal nets and back}
Let $V$ be a unitary VOA. An eigenvector $a$ of $L_0$ is said to be a {\it homogeneous} vector and the corresponding eigenvalue $d$ is called the conformal weight of $a$. If $a \in V$ is homogeneous with conformal weight $d$ one can define the endomorphisms $a_n := a_{(n+d-1)}$, $n\in \mathbb{Z}$, so that $Y(a,z)=\sum_{\mathbb{Z}}a_n z^{-n-d}$. The conformal vector $\nu \in V$ is homogeneous with weight $2$ and $\nu_n=L_n$ for all 
$n\in \mathbb{Z}$. An arbitrary vector $a\in V$ is a finite sum of homogeneous vectors and one defines $a_n$ by linearity. We say that $V$ is {\it energy-bounded} if for every $a \in V$ there exist positive integers $s,k$ and a constant $M>0$ such that 
$\|a_nb\| \leq M(|n| +1)^{s} \|(L_0 + 1_V)^kb\|$ for all $n \in \mathbb{Z}$ and all $b\in V$, where $\|\cdot\|$ is the norm induced on $V$ by the scalar product . 

Now, assume that $V$ is energy-bounded and let $\mathcal{H}_V\supset V$ be its Hilbert space completion.
Then, for every smooth function $f \in C^\infty(S^1)$ and every $a\in V$, we can define the operator $Y_0(a,f)$ on $\mathcal{H}_V$ with domain $V$ by $Y_0(a,f)b=\sum_{n\in \mathbb{Z}}\hat{f}_na_nb$, $b\in V$, where $\hat{f}_n$ denotes the $n$th Fourier coefficient of $f$. Each $Y_0(a,f)$ 
is a closable operator. This means that the closure of its graph is the graph of a possibly unbounded operator $Y(a,f)$ which is by definition a closed operator, i.e.~it has a a closed graph. We call the operators $Y(a,f)$ the {\it smeared vertex operators}.

Given a family $\mathcal{S}$ of (possibly unbounded) closed operators on a Hilbert space $\mathcal{H}$ one can define the von Neumann algebra 
$W^*(\mathcal{S})$ generated by $\mathcal{S}$. A bounded operator $A$ on $\mathcal{H}$ belongs to $W^*(\mathcal{S})$ if and only if $A$ commutes with every unitary operator $U$ such that $UTU^*=T$ for all $T\in \mathcal{S}$.   
Now, for every interval $I\in \mathcal{I}$ we consider the von Neumann algebra 
$$\mathcal{A}_V (I) \equiv W^*(\{Y(a,f): a\in V, f \in C^\infty_c(I) \})$$
generated by all the vertex operators smeared with test functions with support in $I$. 
Since the smeared vertex operators are in general unbounded it is not \textit{a priori} clear that the locality axiom for the VOA $V$ implies that the the map $I \mapsto \mathcal{A}_V(I)$ will satisfy the locality axiom for conformal nets. We say that $V$ is {\it strongly local} if this is actually the case. 
If $V$ is a strongly local VOA then the map $I \mapsto \mathcal{A}_V(I)$ defines an irreducible conformal net $\mathcal{A}_V$ on $S^1$. 

The class of strongly local VOAs turns out to be closed under taking tensor products and unitary subVOAs. Moreover, for every strongly local VOA 
$V$, the map $W\mapsto \mathcal{A}_W$ gives a one-to-one correspondence between the unitary subVOAs $W \subset V$ and the covariant subnets of $\mathcal{A}_V$. Furthermore, if the VOA automorphism group $\mathrm{Aut}(V)$ of a strongly local VOA $V$ is finite then, it coincides with the automorphism group $\mathrm{Aut}(\mathcal{A}_V)$ of the corresponding conformal net $\mathcal{A}_V$. 
Many known examples of unitary VOAs can be proved to be strongly local. This is e.g.~the case for  the unitary affine Lie algebra VOAs, the unitary Virasoro VOAs, 
the known $c=1$ unitary VOAs, the moonshine VOA $V^\natural$, the discrete series $N=1$ and $N=2$ even super-Virasoro VOAs, together with the associated coset and orbifold subVOAs. 
The corresponding conformal nets coincide with those previously constructed by different methods: the loop groups nets \cite{FrG,W},
the Virasoro nets \cite{C2004,KL04}, the $c=1$ conformal nets \cite{Xu2005}, the moonshine conformal net $\mathcal{A}^\natural$ \cite{KL06}, the discrete series $N=1$ and $N=2$ even super-Virasoro nets \cite{CKL,CHKLX}, and the associated coset and orbifold conformal nets \cite{Xu2000,Xu2000b}. 

The even shorter moonshine vertex operator algebra constructed by H\"{o}hn \cite{H} also turns out to be strongly local being a subVOA of 
$V^\natural$. Moreover, the automorphism group of the corresponding conformal net coincides the VOA automorphism group which is known to be the Baby Monster group $\mathbb{B}$, the second largest among the sporadic finite simple groups. This gives for the first time a construction of a conformal net whose automorphism group is $\mathbb{B}$. Furthermore, the (still hypothetical) Haagerup VOA with $c=8$ considered by Evans and Gannon \cite{EG}, related to the Haagerup subfactor, has been suggested to be a unitary subVOA of a unitary affine Lie algebra VOA and hence it should be strongly local. 

To get back a strongly local VOA $V$ from the corresponding conformal net $\mathcal{A}_V$ one can successfully use a construction by Fredenhagen and J\"{o}r{\ss} \cite{FJ}. More generally, it is shown in \cite{CKLW} that the existence of suitable energy bounds for a conformal net $\mathcal{A}$ implies that $\mathcal{A} = \mathcal{A}_V$ for some strongly local vertex operator algebra $V$. 
We conjecture that every unitary VOA is strongly local and that every conformal net comes from a unitary VOA in the way described above. 

\section{Outlook} One of the most interesting aspects of the structural similarities between conformal nets and VOAs is their representation theory. 
In various chiral CFT models the representation theory of the corresponding VOAs and conformal nets can be directly related and have similar properties, e.g.~the same fusion rules. The results described in this contribution appear to be a crucial preliminary step towards the understanding of this relation from a general point of view. We plan to come to the representation theory aspects of the correspondence  $V \mapsto \mathcal{A}_V$ in the near  future. 

\section*{Acknowledgments}
The author is supported in part by the ERC advanced grant 669240 QUEST ``Quantum Algebraic Structures and Models'',
PRIN-MIUR and GNAMPA-INDAM.

\end{document}